\documentclass[letterpaper, 10 pt, journal, twoside]{ieeeconf}	

\usepackage{comment}
\usepackage{graphics} 
\usepackage{epsfig} 
\usepackage{mathptmx} 
\usepackage{times} 
\usepackage{amsmath,bm,mathtools} 
\usepackage{amssymb}  
\usepackage[utf8]{inputenc}
\usepackage[T1]{fontenc}
\usepackage[section]{placeins}
\usepackage{dsfont}
\usepackage{algpseudocode}
\usepackage{algorithm}
\usepackage{url}
\usepackage{tabularx}
\usepackage[percent]{overpic}
\usepackage{cite} 
\usepackage{calrsfs}
\DeclareMathAlphabet{\pazocal}{OMS}{zplm}{m}{n}
\usepackage{color}

\newtheorem{definition}{Definition}
\newtheorem{assumption}{Assumption}

\newtheorem{lemma}{Lemma}
\newtheorem{thrm}{Theorem}
\newtheorem{cor}{Corollary}

\newtheorem{remark}{Remark}

\newcommand\ddfrac[2]{\frac{\displaystyle #1}{\displaystyle #2}}

\usepackage{xcolor,calc}

\begin{document}
%

\title{Distributionally robust stability of payoff allocations in stochastic coalitional games}

\author{
	George Pantazis, Barbara Franci, Sergio Grammatico and Kostas Margellos
	\thanks{The authors are with Delft Center for Systems \& Control, TU Delft  ({\tt \footnotesize \{g.pantazis, s.grammatico\}@tudelft.nl}), the Department of Advanced Computing Sciences in Maastricht University ({\tt \footnotesize b.franci@maastrichtuniversity.nl}) and the Department of Engineering Science, University of Oxford,  ({\tt \footnotesize kostas.margellos@eng.ox.ac.uk}). This work has been partially funded by the ERC under project COSMOS (802348). \\
	©2023 IEEE. Personal use of this material is permitted. Permission from IEEE must be obtained for all other uses, in any current or future media, including reprinting/republishing this material for advertising or promotional purposes, 
	creating new collective works, for resale or redistribution to servers or lists, or reuse of any copyrighted component of this work in other works.
}
}

\maketitle
\thispagestyle{empty}
\pagestyle{empty}

\begin{abstract} 
	
	We consider multi-agent coalitional games with uncertainty in the coalitional values. We provide a novel methodology to study the stability of the grand coalition in the case where each coalition constructs ambiguity sets for the (possibly) unknown probability distribution of the uncertainty. As a less conservative solution concept compared to worst-case approaches for coalitional stability, we consider a stochastic version of the so-called core set, i.e.,  the expected value core. Unfortunately, without exact knowledge of the probability distribution, the evaluation of the expected value core is an extremely challenging task. Hence, we propose the concept of  distributionaly robust (DR) core. 
Leveraging tools from data-driven DR optimization under the Wasserstein distance, we provide finite-sample guarantees that any allocation which lies in the DR core is also stable with respect to the true probability distribution. Furthermore, we show that as the number of samples grows unbounded, the DR core converges almost surely to the true expected value core.  We dedicate the last section to the computational tractability of finding an allocation in the DR core.

\end{abstract}


\section{Introduction}
Coalitional games \cite{Chalkiadakis} are prevalent in  applications ranging from engineering \cite{Raja, Malcolm2, Fele_coalitional} to economics  and social sciences \cite{Mccain2008}. Even though agents in such systems typically act as selfish entities, they are incentivized to form coalitions aiming at receiving higher individual gains or reducing their own  costs. A challenging task, due to the agents' individual interests, is to distribute their payoffs in such a way that none of them has an incentive to deviate from the so-called \emph{grand coalition}, i.e., the coalition where all agents work together. In the literature of coalitional game theory this problem is known as \emph{stability of the grand coalition} and the set of payoffs for which stability is achieved is known as the \emph{core} of the game. Due to its conceptual simplicity, the core has been widely used as a stability concept in coalitional games  \cite{Chalkiadakis} and in turn intense research has been dedicated to finding allocations that lie within the core. \par
Stability of the grand coalition is fundamentally connected to the values of each coalition. 
However,  coalitional values are typically subject to uncertainty. As such, the mathematical framework of deterministic coalitional games needs to be revisited and extended. The seminal works \cite{Charnes_Chance_I,Charnes_Chance_II,Charnes1973} are the first on stochastic coalitional games. The work in \cite{Suijs1} also studies uncertain coalitional games and shows that for a particular class, certain properties of the game, such as the non-emptiness of the core continue to hold when uncertainty is introduced. Uncertain coalitional games were studied under the lenses of Bayesian learning in \cite{Chalkiadakis_RL}, \cite{Ieon_Bayesian}, while the work in \cite{Yuqian2015} investigates which stability solution concepts maximize the probabilistic stability of allocations after the samples of the uncertainty have been revealed.  Moreover, \cite{Repeated2} and \cite{Repeated1} focus on the dynamic evolution of repeated stochastic games. \par  
In \cite{Raja},  the concept of the so called \emph{robust core} is introduced as a generalization to the traditional deterministic core.  In this setting, the range of possible coalitional values is assumed to be known.
The work in \cite{Pantazis_2022_coalitional} extends the notion of the robust core  to that of the \emph{scenario core} accounting for the more general case where both the support set and the probability distribution of the uncertainty affecting the coalitional value are unknown. As an alternative to the robust core \cite{Raja} and its data-driven counterpart \cite{Pantazis_2022_coalitional}, in this paper we consider instead  the significantly milder concept of stability in the mean sense that in turn gives rise to  the so-called \emph{expected value core}. Studying allocation stability in the mean sense might circumvent the possible emptiness of the core set, a fundamental technical challenge in coalitional game theory. Apart from very mild assumptions on the probability distribution of the uncertainty, here we consider both the support set and the probability distribution to be unknown. In other cases, the uncertain parameter affecting the coalitional game might not even admit a single distribution, but a range of possible distributions, quantified through data-driven approaches. As such, evaluating the expected core in this setting is extremely challenging. \par 
To address this, we follow an approach based on distributionally robust (DR) optimization \cite{shapiro_minimax_2002, fournier_rate_2015, calafiore_distributionally_2006, goh_distributionally_2010, delage_distributionally_2010, cherukuri_cooperative_2020, gao_distributionally_2022} thus considering ambiguity sets that represent empirical sets in which the true probability distribution (in case the uncertainty admits one) is likely to be contained. The consideration of ambiguity sets leads to allocations that are \emph{distributionally stable}. We call the set of all distributionally stable allocations the \emph{distributionally robust (DR)} core of the game. \par 
	   Leveraging results from data-driven DR optimization under the Wasserstein distance \cite{mohajerin_esfahani_data-driven_2018},  \cite{netessine_wasserstein_2019}, we provide finite sample guarantees on the probability that any allocation in the DR core of the DR game approximation is also in the expected value core  of the original game with a given confidence (Section III). 
	      Moreover, we prove almost-sure asymptotic convergence  of the Wasserstein DR core to the expected core of the original game (Section IV.A). Finally, we provide the means to calculate an allocation in the Wasserstein DR core (Section IV.B). Specifically, we show that under certain conditions, the problem of finding  such an allocation can be recast as a convex optimization problem, whose complexity both in the number of decision variables and constraints is inherently connected to the number of possible subcoalitions. Numerical simulations corroborate our theoretical findings (Section V).    \par 
\section{Stochastic coalitional games}
\subsection{Allocation mean stability and expected value core }

We consider a coalitional game with $N$ agents parameterized by the index set $\pazocal{N}=\{1, \dots, N\}$.  We denote the number of possible subcoalitions except for the grand coalition by $M$, i.e., $M=2^N-1$. In this setting, the agents, though selfish, wish to form coalitions if that implies an increase in their individual payoffs. The total gain for each coalition is given by the so called \emph{value function}, which, depending on the coalition $S \subseteq \pazocal{N}$, takes a real value  representing the total payoff that agents participating in it would obtain from its formation.  However, the values of each coalition are subject to uncertainty thus rendering the value function  of each coalition stochastic. 
\begin{definition} \emph{(Stochastic value function)}  
	The value function of a coalition $S \subset \pazocal{N}$ is a function $u_S: 2^\pazocal{N} \times \Xi \rightarrow \mathbb{R}$ that, given the value of the uncertainty realization  $\xi \in \Xi \subseteq \mathbb{R}^p$, returns  the total payoff for the agents forming a coalition $S$.  The value function of the grand coalition is deterministic, i.e., $u_\pazocal{N}: 2^\pazocal{N} \rightarrow \mathbb{R}$. \hfill $\square$
	\end{definition}
An uncertain coalitional game is then defined as the tuple $G_\mathbb{P}=\{\pazocal{N},\{u_S\}_{S \subseteq \pazocal{N}}, \Xi, \mathbb{P}\}$, where $\mathbb{P}$ denotes the probability distribution that the uncertainty $\xi \in \Xi$ follows. 
To circumvent the fundamental issue of the emptiness of the robust core as defined in \cite{Raja}, let us consider  the concept of stability of allocations in the mean sense, defined as follows.
\begin{definition} \label{mean_allocation} \emph{(Stability in the mean sense)} 
	An allocation  $x=(x_i)_{i \in \pazocal{N}}$ of the game $G_\mathbb{P}=\{\pazocal{N},\{u_S\}_{S \subseteq \pazocal{N}}, \Xi, \mathbb{P}\}$ is \emph{stable in the mean sense }  if 
i) $\sum_{i \in \pazocal{N}} x_i = u_\pazocal{N}$ and 
	ii) $\sum_{i \in S}x_i \geq \mathbb{E}_{\mathbb{P}}[u_S(\xi)],  \ \forall \ S \subset \pazocal{N}$.  \hfill $\square$
\end{definition}
The first condition is called the efficiency condition. Due to our assumption that the grand coalition is deterministic, it means that the total increase in gains when all agents work together is known with certainty.  
This is the case when agents might know how efficient a fully-cooperative scheme is but have some level of uncertainty/ambiguity with respect to the potential outcomes of the subcoalitions. The second condition implies that the allocation $x$  is not strictly feasible, hence agents do not have an incentive to form $S$. Otherwise, if $\sum_{i \in S}x_i < \mathbb{E}_{\mathbb{P}}[u_S(\xi)]$  agents would have the incentive to leave the grand coalition and form $S$, thus receiving $ \mathbb{E}_{\mathbb{P}}[u_S(\xi)]$ as opposed to $\sum_{i \in S}x_i$. In this setting, we wish to study the stability of the grand coalition, where no agent has an incentive to deviate and create other subcoalitions. To this end, let us introduce an extension to the classic notion of the  core, the so-called \emph{expected value core},  which is the set of all stable allocations in the mean sense as defined next. 
\begin{definition} \label{mean_core} \emph{(Expected value core)}
	The expected value core $C_E(G_{\mathbb{P}})$ of the game  $G_{\mathbb{P}}$ is defined as  the set 
	\begin{align}
		\qquad C_{E}(G_{\mathbb{P}})= \{ x \in \mathbb{R}^N:& \sum\limits_{i \in \pazocal{N}} x_i = u_\pazocal{N}, \nonumber \\ 
		& \sum\limits_{i \in S}x_i \geq \mathbb{E}_{\mathbb{P}}[u_S(\xi)],  \ \forall \ S \subset \pazocal{N}\},\qquad  \nonumber 
	\end{align} 
where
\begin{align}   \label{expected_integral}
	\mathbb{E}_\mathbb{P}[u_S(\xi)]= \int_{\Xi}u_S(\xi)\mathbb{P}(d\xi).  
\end{align}
\end{definition} \hfill $ \square$ \\
 We now impose the following mild technical assumptions.
 \begin{assumption} \emph{(Light-tailed probability distribution)} \label{light_tail}
 	For the true probability measure $\mathbb{P}$, there exists $a >1$ such that 
 	\begin{align}
 		A:=		\mathbb{E}_\mathbb{P}[\exp(\|\xi\|^a)]= \int_{\Xi}\exp(\|\xi\|^a) \mathbb{P}(d\xi) < \infty. \nonumber 
 	\end{align}
 \end{assumption}
\begin{assumption} \label{finite_moment}
	For the probability measure $\mathbb{P}$ it holds that $\mathbb{E}_\mathbb{P}[ \|\xi\| ]=\int_\Xi \|\xi \| \mathbb{P}(d\xi) < \infty$.
\end{assumption}
 Assumption \ref{light_tail} requires the tail of the true probability distribution $\mathbb{P}$ to decay at an exponential rate. In case $\Xi$ is a compact set this assumption trivially holds. Assumption \ref{finite_moment} requires that $\mathbb{P}$ admits a finite first-order moment. 
 In such a general set-up it is challenging, if not impossible, to compute the expected value integral given in (\ref{expected_integral}), i.e., one cannot evaluate the expected-valued core $C_{E}(G_{\mathbb{P}})$.
 To circumvent this challenge, we propose a methodology based on distributionally robust optimization. 

\subsection{Distributionally robust stability of allocations}
  In our setting, agent coalitions  $S \subset \pazocal{N}$ can construct ambiguity sets of the probability distribution $\mathbb{P}$ of the uncertainty $\xi \in \Xi$ that affect their coalitional values $u_S(\xi)$. This is due to lack of knowledge of  $\mathbb{P}$. In other words, we do not only have uncertainty affecting the coalitional game, but uncertainty about the distribution of the uncertain parameter.   \par 
 We postulate that each coalition $S \subset \pazocal{N}$ is allowed to construct their own ambiguity sets. The heterogeneity of the coalitional ambiguity sets provides the necessary modelling freedom for our theory to be flexible for application purposes. 
To this end, we assume that each coalition $S \subset \pazocal{N}$ has access to their own i.i.d. samples $\xi_{K_S} = (\xi^{(1)}, \dots, \xi^{(K_S)}) \in \Xi^{K_S}$ and consider the distributionally robust version $G_{\pazocal{\hat{P}}_K}$ of the original game $G_{\mathbb{P}}$ defined as the tuple  $G_{\pazocal{\hat{P}}_K}=\{\pazocal{N},\{u_S\}_{S \subseteq \pazocal{N}}, \Xi, \pazocal{\hat{P}}_K\}$, where $K= \{K_S\}_{S \subset \pazocal{N}}$, while  $\pazocal{\hat{P}}_K= \{\pazocal{\hat{P}}_{K_S}\}_{S \subset \pazocal{N}}$ is the collection of ambiguity sets constructed based on the available data $\xi_{K_S}$ of each subcoalition $S \subset \pazocal{N}$. We now proceed to defining the notion of distributional stability of an allocation. 
\begin{definition} \emph{(Distributionally robust  stability of allocations)} 
	For a given number of i.i.d. drawn samples $ \xi_{K_S}=(\xi^{(1)}, \dots, \xi^{(K_S)}) \in \Xi^{K_S}$ per coalition $S \subset \pazocal{N}$, an allocation $x=(x_i)_{i \in \pazocal{N}}$ is distributionally stable with respect to the coalitional ambiguity sets $\pazocal{\hat{P}}_{K_S}$, $S \subset \pazocal{N}$ if
	\begin{enumerate}
		\item $\sum_{i \in \pazocal{N}} x_i = u_\pazocal{N}$ and 
		\item $\sum_{i \in S}x_i \geq \sup_{\mathbb{Q}_S \in \pazocal{\hat{P}}_{K_S}}\mathbb{E}_{\mathbb{Q}}[u_S(\xi)],  \ \forall \ S \subset \pazocal{N}$.  \hfill $\square$
	\end{enumerate}
\end{definition}

\section{Distributionally robust coalitional games based on the Wasserstein distance}
\subsection{Background on distributional robustness}
In this section, we introduce some basic concepts from distributionally robust optimization under the Wasserstein metric \cite{fournier_rate_2015}, \cite{mohajerin_esfahani_data-driven_2018}, \cite{netessine_wasserstein_2019}.
We show how one can leverage this framework in order to provide certificates of stability with respect to the true unknown coalitional game along with a tractable approximation of its expected value core. We start by  imposing a mild assumption on the probability distributions serving as candidates of the true distribution $\mathbb{P}$. Specifically, we consider all distributions with bounded first-order moments, i.e., $\mathbb{Q} \in \pazocal{M}(\Xi)$, where  $\pazocal{M}(\Xi)$ is the set of probability distributions with support $\Xi$ that satisfy Assumption \ref{finite_moment}.
We then need a measure of distance between two probability distributions to quantify how close a candidate probability distribution is to the true probability distribution $\mathbb{P}$; Let us thus use the Wasserstein distance defined as follows. 

\begin{definition} \label{Wasserstein}(\emph{Wasserstein distance}, \cite{Kantorovich_1958})
	The Wasserstein distance $d_W: \pazocal{M}(\Xi) \times \pazocal{M}(\Xi) \rightarrow \mathbb{R}_{\geq 0}$  between two probability distributions $\mathbb{Q}_1, \mathbb{Q}_2 \in \pazocal{M}(\Xi)$ is defined as 
	\begin{align}
		d_W(\mathbb{Q}_1, \mathbb{Q}_2)=&\inf_{\Pi}\Big\{\int_{\Xi^2} \|\xi_1-\xi_2 \| \Pi(d\xi_1, d\xi_2): \nonumber \\ 
		& \Pi \text{ is a joint distribution of $\xi_1$ and $\xi_2$}  \nonumber \\
		&\text{	with marginals $\mathbb{Q}_1$ and $\mathbb{Q}_2$, respectively}\Big\}, \nonumber 
	\end{align}
	 where $\| \cdot \|$ can be any norm on $\mathbb{R}^p$. 
\end{definition}
An alternative dual interpretation can be derived by the so-called Kantorovich-Rubinstein  theorem:
\begin{thrm}   \label{KR_theorem} (\emph{Kantorovich-Rubinstein theorem}, \cite{Kantorovich_1958})
	Given two probability distributions $Q_1$, $Q_2\in \pazocal{M}(\Xi)$ we have that 
	\begin{align}
	d_W(\mathbb{Q}_1, \mathbb{Q}_2)= \sup\limits_{f \in \pazocal{F}}\left\{ \int_{\Xi} f(\xi) \mathbb{Q}_1(d\xi)- \int_{\Xi} f(\xi) \mathbb{Q}_2(d\xi) \right\},  \nonumber 
		\end{align}
	where $\pazocal{F}$ is the space of all Lipschitz continuous functions for which $|f(\xi)-f(\xi')| \leq ||\xi- \xi'||$ for all $\xi, \xi' \in \Xi$. \hfill $\square$
	\end{thrm}
\subsection{ Finite-sample guarantees for the distributionally robust core}

In the subsequent developments, we consider
 that each coalition $S \subset \pazocal{N}$ has their own independent samples from $\mathbb{P}$. Any given coalition $S \subset \pazocal{N}$  constructs their respective ambiguity set based on their collected data $\xi_{K_S}=(\xi^{(k_S)})_{k_S=1}^{K_S} \in \Xi^{K_S}$. For each coalition $S \subset \pazocal{N}$ each ambiguity set is given by:
 \begin{align}
 	\mathbb{B}_{\epsilon_S}(\mathbb{\hat{P}}_{K_S})=\{\mathbb{Q}_S \in \pazocal{M}: d_W(\mathbb{\hat{P}}_{K_S}, \mathbb{Q}_S) \leq \epsilon_S\},
 \end{align}
where  $\mathbb{\hat{P}}_{K_S}=\sum_{k_s=1}^{K_S} \delta(\xi-\xi^{(k_S)})$  is the empirical probability distribution of each coalition $S$  on the basis of $K_S$ i.i.d samples  from the support set $\Xi$ by coalition $S$ with $\delta(\xi-\xi^{(k_S)})=\frac{1}{K_S}$ if  $\xi = \xi^{(k_S)}$, $k_S=\{1, \dots, K_S\}$ and $0$ otherwise. 

The following result relies on Assumption \ref{light_tail} to provide guarantees on the probability that a multi-sample will be drawn from coalition $S \subset \pazocal{N}$ such that true probability measure lies within the constructed Wasserstein ball with a given confidence.

\begin{lemma} \label{DR_guarantees} 
Let Assumption \ref{light_tail} hold and for any coalition $S \subset \pazocal{N}$ fix $\epsilon_S>0$. We have that 
	\begin{align}
		\mathbb{P}^{K_S}\left\{\xi_{K_S} \in \Xi^{K_S}: d_W(\mathbb{P}, \mathbb{\hat{P}}_{K_S})\leq \epsilon_S \right\} \geq 1-\beta_S, \nonumber 
	\end{align}
	where 
	\begin{align} \label{beta}
		\beta_S= \begin{cases} 
			c \exp(-qK_S\epsilon_S^{\max{\{p,2\}}}), & \text{ if } \epsilon_S \leq 1 \\
			c \exp(-q K_S\epsilon_S^{a}),  &  \text{ if } \epsilon_S >1,  
		\end{cases}
	\end{align}
	for all $K _S  \geq 1$ and $p \neq 2$, where $c, q$ are positive constants that only depend on  the parameters $a, A$ in Assumption \ref{light_tail} and the dimension of the support set $p$.
\end{lemma}
\emph{Proof}: The proof is an adaptation of Theorem 2 in \cite{fournier_rate_2015} applied to each Wasserstein ball  $\pazocal{\hat{P}}_{K_S}=\mathbb{B}_{\epsilon_S}(\mathbb{\hat{P}}_{K_S})$ constructed around the empirical probability distribution $\mathbb{\hat{P}}_{K_S}$ and a radius $\epsilon_S$ chosen by each coalition $S \subset \pazocal{N}$. \hfill $\blacksquare$
Lemma \ref{DR_guarantees} paves the way towards establishing finite sample guarantees for the following Wasserstein-based version of the distributionally robust core.
\begin{definition} \emph{(Distributionally robust core)}\label{DR_core}
	The distributionally robust core $C_{\text{DR}}( G_{\pazocal{\hat{P}}_K} )$ of the game  $G_{\mathbb{P}}$ based on the Wasserstein distance is defined as  the set 
	\begin{align}
		C_{\text{DR}}( G_{\pazocal{\hat{P}}_K} ) &= \Big\{ x \in \mathbb{R}^N: \sum\limits_{i \in \pazocal{N}} x_i = u_\pazocal{N}, \nonumber \\ 
		& \sum_{i \in S}x_i \geq \sup_{Q_S \in \pazocal{\hat{P}}_{K_S}}\mathbb{E}_{Q_S}[u_S(\xi_S)],  \ \forall \ S \subset \pazocal{N} \Big\},  \nonumber 
	\end{align}
	where $\pazocal{\hat{P}}_{K_S}=\mathbb{B}_{\epsilon_S}(\mathbb{\hat{P}}_{K_S})$. \hfill $\square$
\end{definition}

Throughout we assume that for all multi-samples, the ambiguity sets are such that a non-empty DR core is returned.

\begin{thrm} \label{main_theorem_dr}
	For each $S \subset \pazocal{N}$ fix a Wasserstein radius $\epsilon_S$ and consider a multi-sample size $K_S$. It holds that 
	\begin{align}
		\mathbb{P}^K \Big\{ \xi_{K} \in \Xi^K: C_{\text{E}}( G_\mathbb{P} ) \supseteq &C_{\text{DR}}( G_{\pazocal{\hat{P}}_K} ) \Big\} \geq \beta , \label{prob_guarantess}
	\end{align}
where $\beta= \prod_{S \subset \pazocal{N}}(1- \beta_S)$ and each $\beta_S$ is given by (\ref{beta}).
\end{thrm}
 \emph{Proof}: 
  We have that 
\begin{align}
	&\mathbb{P}^K\Big\{\xi_K \in \Xi^K: \ C_{\text{E}}( G_{\{\Xi, \mathbb{P}\}} ) \supseteq C_{\text{DR}}( G_{\{\Xi, \mathbb{P}\}})\Big\}  \nonumber  \\
	&\geq \mathbb{P}^K \Big\{\xi_K \in \Xi^K:  \mathbb{E}_{\mathbb{P}} [u_S(\xi)] \leq \sup_{\mathbb{Q}_S \in \pazocal{\hat{P}}_{K_S}} \mathbb{E}_{\mathbb{Q}_S} [u_S(\xi)], \ \forall \ S \subset \pazocal{N} \Big\} \nonumber \\
	&= \mathbb{P}^K\Big\{ \bigcap_{S \subset \pazocal{N}}\Big\{\xi_K \in \Xi^K:  \mathbb{E}_{\mathbb{P}} [u_S(\xi)] \leq \! \! \!  \! \! \!  \sup_{\mathbb{Q}_S \in \pazocal{\hat{P}}_{K_S}} \! \! \!  \! \! \! \mathbb{E}_{\mathbb{Q}_S}  [u_S(\xi)] \Big\} \Big\}   \nonumber \\
	 & = \prod_{S \subset \pazocal{N}}\mathbb{P}^{K_S} \Big\{ \xi_{K_S} \in \Xi^{K_S}:  \mathbb{E}_{\mathbb{P}} [u_S(\xi)] \leq \! \! \! \! \! \sup_{\mathbb{Q}_S \in \pazocal{\hat{P}}_{K_S}} \! \! \! \! \! \! \mathbb{E}_{\mathbb{Q}}  [u_S(\xi)] \Big\}  \label{relation_5}
\end{align}
 The last equality is due to the fact that each coalition constructs its ambiguity sets based on its own (independent) samples.  From Lemma  \ref{DR_guarantees} for each coalition $S \subset \pazocal{N}$ we have that 
\begin{align}
	\mathbb{P}^{K_S}\Big\{\xi_{K_S} \in \Xi^{K_S}: d_W(\mathbb{P}, \mathbb{\hat{P}}_{K_S})\} \leq \epsilon_S \Big\} \geq 1-\beta_S. \nonumber 
\end{align}
Therefore, 
\begin{align}
	&\prod_{S \subset \pazocal{N}}\mathbb{P}^{K_S}\Big\{ \xi_{K_S} \in \Xi^{K_S}:  \mathbb{E}_{\mathbb{P}} [u_S(\xi)] \leq  \! \! \! \! \!  \sup_{\mathbb{Q}_S \in \pazocal{\hat{P}}_{K_S}}   \! \!  \! \! \! \mathbb{E}_{\mathbb{Q}_S}  [u_S(\xi_S)] \Big\}     \nonumber \\  
	&\geq  \prod_{S \subset \pazocal{N}} \mathbb{P}^{K_S}\Big\{\xi_{K_S} \in \Xi^{K_S}: d_W(\mathbb{P}, \mathbb{\hat{P}}_{K_S}) \leq \epsilon_S \Big\}  \nonumber  \\ 
	&\geq \prod_{S \subset \pazocal{N}}(1-\beta_S),  \label{prob_guarantess_agent_type}
\end{align}
where the first inequality follows from Theorem 3.5 in \cite{mohajerin_esfahani_data-driven_2018}.  \hfill $\blacksquare$ 

 The following result provides guarantees when all coalitions use the same parameters.
\begin{cor} \label{homogeneous_guarantess}
Consider Assumption \ref{light_tail}. For each coalition $S \subset \pazocal{N}$ fix a common Wasserstein radius $\epsilon$ and assume that the same multi-sample $\xi_K$ is used among coalitions. Then, it holds that
	\begin{align}
		\mathbb{P}^K\{\xi_K \in \Xi^K:  \ C_{\text{E}}( G_\mathbb{P} )  \supseteq C_{\text{DR}}( G_{\pazocal{\hat{P}}_{K}})\} 	\geq (1 - \beta)^M, \nonumber 
	\end{align}
	where $\beta$ is given by (\ref{beta}) by setting $\epsilon_S=\epsilon$ and $K_S=K$.  \hfill $\square$
\end{cor}
\emph{Proof}: Since $\epsilon_S=\epsilon$ for all $S \subset \pazocal{N}$ and the same number of samples $K \in \mathbb{N}$ is drawn for all $S \subset \pazocal{N}$, then by  (\ref{beta}), $\beta_S=\beta$ for all $S \subset \pazocal{N}$. By  Lemma \ref{DR_guarantees} we then have that for each coalition $S \subset \pazocal{N}$
	\begin{align}
	\mathbb{P}^K\Big\{\xi_K \in \Xi^K:  d_W(\mathbb{P}, \mathbb{\hat{P}}_{ K}) \leq \epsilon\Big\} 	\geq 1 - \beta, \nonumber 
\end{align}
where $\beta$ is given by (\ref{beta}).
 Following the same proofline as in Theorem \ref{main_theorem_dr} and setting $\beta_S=\beta$ to the right hand side of (\ref{prob_guarantess_agent_type}) concludes the proof. \hfill $\blacksquare$ \par 
  \begin{remark}
 	The confidence parameter in Theorem \ref{main_theorem_dr} and Corollary 1  depends on the number of possible subcoalitions $S \subset \pazocal{N}$. Corollary 1 shows that the dependence on the number of possible coalitions of the original problem is also inherited by the provided finite-sample guarantees. Though under a different approach, this is also observed in the \emph{a priori} results in \cite{Pantazis_2022_coalitional}, where the authors apply the results of  \cite{pantazis_posteriori_2020}, \cite{Pantazis2020} to construct a data-driven version of the robust core.   \hfill $\square$
 \end{remark}
 \subsection{Agent-level sampling}
  Assume now that the ambiguity sets are constructed on the basis of samples drawn by the individual agents, i.e, with a slight abuse of notation each agent draws a multi-sample $\xi_{K_i}$, which is then used by any coalition $S \subset\pazocal{N}$ for which $i \in S$.  Fix a confidence parameter $\bar{\beta}_S \in (0,1)$ for each $S \subset \pazocal{N}$ and let $\xi_{ \bar{K} }= (\xi_{\bar{K}_S})_{S \subset \pazocal{N}}$, where $\xi_{\bar{K}_S}$ are the samples of each coalition $S \subset \pazocal{N}$. Then, it is easy to show that 
 \begin{align}
 	\mathbb{P}^{\bar{K}}\{\xi_{\bar{K}} \in \Xi^{\bar{K}}: C_{\text{E}}( G_\mathbb{P} ) \supseteq &C_{\text{DR}}( G_{\pazocal{\hat{P}}_{\bar{K}}} )\} \geq  \bar{\beta}, \nonumber 
 \end{align}
 where $\bar{\beta}=\sum_{S \subset \pazocal{N}}(1- \bar{\beta}_S)-M+1$ and each $\bar{\beta}_S$ is given by (\ref{beta}) setting  $\epsilon_S=\bar{\epsilon}_S$ and $K_S=\bar{K}_S= \sum_{i \in S}K_i$.  This result is obtained by applying Bonferroni's inequality to the third step of the proof of Theorem \ref{main_theorem_dr} and thus constitutes a rather conservative bound. Further work is required to leverage the data on the agents' level and translate it to guarantees on a coalitional level taking into account sharing of data among coalitions and thus improving those theoretical guarantees. This is a challenging task and thus left for future work.

\section{Asymptotic consistency and computational tractability of the DR core}
\subsection{Asymptotic consistency of the DR core}
In this subsection we show that under Lipschitz continuity of the value functions and careful choice of the radius of the Wasserstein ball  $\epsilon_S$ and of the confidence parameter $\beta_S$ for each coalition $S \subset \pazocal{N}$, the DR core based on the Wasserstein distance converges almost surely to the true expected value core of the original problem. Let us impose the following assumption: 
\begin{assumption} \label{lipschitz}
	For each coalition $S \subset \pazocal{N}$ the value function $u_S$ is $L_S$-Lipschitz continuous in $\xi$ with $L_S \geq 0$, i.e, $\|u_S(\xi)-u_S(\xi')\| \leq L_S\|\xi-\xi'\|$ for all $\xi, \xi' \in \Xi$.
	 \hfill $\square$
\end{assumption}
Such assumptions are common in stability analysis for stochastic programming \cite{netessine_wasserstein_2019}, \cite{romisch_stability_1991} and also provide interesting insights in our setting.  In fact, the increase in the error between a DR coalitional value and the corresponding  expected coalitional value is proportional to the estimation error between the true and the empirical distribution of that coalition amplified by at most the Lipschitz constant of the corresponding value function. \par 
As opposed to the developments of the previous section, where  we fix the radius $\epsilon_S$ and the number of samples $K_S$ for any $S \subset \pazocal{N}$ and calculate $\beta_S$ based on (\ref{beta}), we now solve (\ref{beta}) with respect to $\epsilon_S$, thus obtaining the Wasserstein radius as a function of the confidence parameter $\beta_S$ and the number of samples $K_S$. In particular, we have that
	\begin{align} 
		\epsilon_S(\beta_S, K_S) =  \begin{cases} 
			\left (\dfrac{\ln(\frac{c}{\beta_S})}{qK_S} \right )^{\frac{1}{ \max\{p,2\} }} & \text{ if } K_S \geq \dfrac{\ln(\frac{c}{\beta_S})}{q}  \\ 
			\left (\dfrac{\ln(\frac{c}{\beta_S})}{qK_S} \right )^{\frac{1}{a}}   &  \text{ if } K_S < \dfrac{\ln(\frac{c}{\beta_S})}{q}. \nonumber 
		\end{cases}
	\end{align}

  The following theorem establishes almost-sure convergence of the DR core to the expected value-core as the number of samples increases.
\begin{thrm} \label{asymptotic_consistency}
	Let Assumptions \ref{light_tail} and \ref{lipschitz} hold. Suppose that for each $S \subset \pazocal{N}$, $\beta^{K_S}_S \in (0,1), K_S \in \mathbb{N}$ satisfies $\sum\limits _{K_S =1}^\infty \beta^{K_S}_{S} < \infty $ and $\lim\limits_{K_S \rightarrow \infty}\epsilon_S(\beta_S, K_S)=0$. Any sequence of the Wasserstein-based  DR cores $\{C_{\text{DR}}( G_{\pazocal{\hat{P}}_{K}})\}_{K \in \mathbb{N}^{M}}$, where $K=(K_S)_{S \subset \pazocal{N}}$, converges $\mathbb{P}^\infty$-almost surely to the true expected core $C_{\text{E}}( G_\mathbb{P} )$ as $K_S \rightarrow \infty$ for all $S \subset \pazocal{N}$. \hfill $\square$ 
\end{thrm}
\emph{Proof}:
For each coalition $S \subset \pazocal{N}$ we have that 
\begin{align}
&	\mathbb{P}^{K_S}\{ \xi_{K_S} \in \Xi^{K_S}:  \mathbb{E}_{\mathbb{P}} [u_S(\xi)] \leq \! \! \! \! \! \! \sup_{\mathbb{Q}_S \in \pazocal{\hat{P}}_{K_S}} \! \! \! \! \! \! \mathbb{E}_{\mathbb{Q}_S}  [u_S(\xi)]\}  \nonumber  \\
	&\geq  \mathbb{P}^{K_S}\{\xi_{K_S} \in \Xi^{K_S}: d_W(\mathbb{P}, \mathbb{\hat{P}}_{K_S})\leq \epsilon_S(\beta^{K_S}_S, K_S)\} \geq 1-\beta^{K_S}_S, \nonumber 
\end{align}
where the last inequality is due to Lemma \ref{DR_guarantees}. 
Letting $K_S \rightarrow \infty$ and since $\lim\limits_{K_S \rightarrow \infty} \beta_S^{K_S}=0$ we have that 
\begin{align}
\mathbb{E}_{\mathbb{P}} [u_S(\xi)] \leq \lim\limits_{K_S \rightarrow \infty} \sup_{\mathbb{Q}_S \in \pazocal{\hat{P}}_{K_S}} \mathbb{E}_{\mathbb{Q}_S}  [u_S(\xi)] \label{relation_3}
\end{align}
$\mathbb{P}^\infty$-almost surely.
Following a methodology similar in spirit to \cite{mohajerin_esfahani_data-driven_2018}, for each coalition $S \subset \pazocal{N}$ and for every $K_S \in \mathbb{N}$ by the definition of supremum for any $\delta_S>0$ there exists $\mathbb{\hat{Q}}_{K_S}$ such that
\begin{align}
	\sup_{\mathbb{Q}_S \in \pazocal{\hat{P}}_{K_S}} \mathbb{E}_{\mathbb{Q}_S}  [u_S(\xi)] \leq \mathbb{E}_{\mathbb{\hat{Q}}_{K_S}}[u_S(\xi)] + \delta_S  \nonumber 
	\end{align}
By the Kantorovich-Rubinstein theorem (Theorem \ref{KR_theorem}), we have that 
\begin{align}
L_S d_W(\mathbb{\hat{Q}}_{K_S}, \mathbb{P}) 
	\geq L_S \left ( \mathbb{E}_{\mathbb{\hat{Q}}_{K_S}}\left[\frac{1}{L_S}u_S(\xi)\right] -  \mathbb{E}_{\mathbb{P}}\left[\frac{1}{L_S}u_S(\xi)\right] \right ), \nonumber 
\end{align} 
since, under Assumption \ref{lipschitz}, $\ddfrac{u_S(\xi)}{L_S}$ is a Lipschitz continuous function with Lipschitz constant less than or equal to 1. The relation above can then be written as 
\begin{align}
 \mathbb{E}_{\mathbb{\hat{Q}}_{K_S}}[u_S(\xi)] \leq \mathbb{E}_{\mathbb{P}}[u_S(\xi)]+ L_S  d_W(\mathbb{\hat{Q}}_{K_S}, \mathbb{P}). \nonumber 
 \end{align}
We then have that:
\begin{align}
	&\lim_{K_S \rightarrow \infty} \sup_{\mathbb{Q}_S \in \pazocal{\hat{P}}_{K_S}} \mathbb{E}_{{\mathbb{Q}}_S}[u_S(\xi)]  \leq  	\lim_{K_S \rightarrow \infty} \mathbb{E}_{\mathbb{\hat{Q}}_{K_S}}[u_S(\xi)]+\delta_S  \nonumber  \\
	& \leq   \lim_{K_S \rightarrow \infty} \left\{\mathbb{E}_{\mathbb{P}}[u_S(\xi)]+L_S d_W(\mathbb{\hat{Q}}_{K_S}, \mathbb{P}) \right\} +\delta_S= \mathbb{E}_{\mathbb{P}}[u_S(\xi)]+\delta_S  \nonumber 
\end{align}
$\mathbb{P}^\infty$-almost surely, since by adapting \cite[Lemma 3.7]{mohajerin_esfahani_data-driven_2018} in our setting, we have that for each $S \subset \pazocal{N}$
\begin{align}
	\lim_{K_S \rightarrow \infty} d_W(\mathbb{P}, \mathbb{\hat{Q}}_{K_S})=0, \  \mathbb{P}^\infty \text{-almost surely}. \nonumber 
\end{align}
Letting $\delta_S \downarrow 0$ we have that 
\begin{align}
\lim_{K_S \rightarrow \infty} \sup_{\mathbb{Q}_S \in \pazocal{\hat{P}}_{K_S}} \mathbb{E}_{{Q}_S}[u_S(\xi)] \leq  \mathbb{E}_{\mathbb{P}}[u_S(\xi)]. \label{relation_4}
\end{align}
From relations (\ref{relation_3})  and (\ref{relation_4}) we have that for any $S \subset \pazocal{N}$
\begin{align}
\lim_{K_S \rightarrow \infty} \sup_{\mathbb{Q}_S \in \pazocal{\hat{P}}_{K_S}} \mathbb{E}_{{Q}_S}[u_S(\xi)] = \mathbb{E}_{\mathbb{P}}[u_S(\xi)], \  \mathbb{P}^\infty \text{-almost surely}.	\nonumber 
\end{align}
As such,  by Definitions \ref{mean_core} and \ref{DR_core}, any sequence of DR cores $\{C_{\text{DR}}( G_{\pazocal{\hat{P}}_{K}})\}_{K \in \mathbb{N}^{M}}$, for which $S \subset \pazocal{N}$, $\beta^{K_S}_S \in (0,1), K_S \in \mathbb{N}$ satisfies $\sum\limits _{K_S =1}^\infty \beta^{K_S}_{S} < \infty $ and $\lim\limits_{K_S \rightarrow \infty}\epsilon_S(\beta_S)=0$ for all $S \subset \pazocal{N}$ with $K=(K_S)_{S \subset \pazocal{N}}$, converges $\mathbb{P}^\infty$-almost surely to the true expected core $C_{\text{E}}( G_\mathbb{P} )$ as $K_S \rightarrow \infty$ for all $S \subset \pazocal{N}$. \hfill $\blacksquare$
\subsection{Finding allocations inside the DR core}
The results of this subsection hold irrespective of what type of sampling we perform. Leveraging results from \cite{mohajerin_esfahani_data-driven_2018}, we show that an allocation inside the DR core can be computed by solving a finite-dimensional convex optimization problem. Here we impose the following assumption.
\begin{assumption} \label{properties}
\begin{enumerate}
	\item For any $S \subset \pazocal{N}$ the value function $u_s(\xi)$ can be written as $u_S(\xi)=\max\limits_{m_S=1, \dots, M_S}u_{m_S}(\xi)$, where $-u_{m_S}(\xi)$ is proper, convex and lower semi-continuous for all $m_S \in \{1, \dots, M_S\}$ and any $S \subset \pazocal{N}$.
	\item For any $S \subset \pazocal{N}$ $u_{S}$ does not take the value $-\infty$  on $\Xi$.
	\item The support set $\Xi$ is closed and convex.  \hfill $\square$
	\end{enumerate}
\end{assumption}
Under these assumptions we have the following result:
\begin{lemma} \label{convex_program}
	Let Assumption \ref{properties} hold.  By drawing $K_S$ samples and considering the dual variables $\lambda_S, \ell_{k_S}, z_{k_Sm_S}, v_{k_Sm_S}$ that correspond to the Wasserstein ball constraint of each coalition $S \subset \pazocal{N}$, an allocation inside the DR core is found by solving the optimization problem
\begin{equation}
P:	\left\{
	\begin{aligned}
	 &	\min_{ \substack{x, \{\lambda_S, \ell_{k_S}, z_{k_Sm_S}, v_{k_Sm_S}\}_{S \subset \pazocal{N}} }}	||x||^2_2 \\
	&	\text{s.t. } \sum_{i \in \pazocal{N}}x_i = u_\pazocal{N} \nonumber  \\
		&\lambda_S \epsilon_S + \frac{1}{K_S}\sum_{k_S=1}^{K_S}\ell_{k_S} \leq \sum\limits_{i \in S}x_i, \  \forall \ S \subset \pazocal{N}   \nonumber  \\
		&   [-u_{m_S}]^\ast(z_{k_Sm_S}-v_{k_Sm_S})+ \sigma_\Xi(v_{k_Sm_S})-z_{k_S}^\top \xi^{(k_S)} \leq \ell_{k_S}, \nonumber \\ 
		& \hspace{5cm} \forall k_S, \forall  m_S, \  \forall \ S  \subset \pazocal{N}   \nonumber \\
		& \  \ \  \   \|z_{k_Sm_S} \|_* \leq \lambda_S \ \forall k_S, \forall  m_S, \  \forall \ S \subset \pazocal{N}  \nonumber 
	\end{aligned}
	\right.
\end{equation} 
where $[f]^*$ denotes the conjugate function of a function $f$, i.e., $[f]^*(y)=\sup_{x \in dom(f)}(y^\top x-f(x))$. and  $\| \cdot \|_{*}$ is the dual norm, while $\sigma_X$ denotes the conjugate of the characteristic function. \hfill $\square$
\end{lemma}
\emph{Proof}: We wish to solve the following feasibility problem
\begin{equation}
	\left\{
	\begin{aligned}
		&\min_{x \in \mathbb{R}^N}	\ ||x||^2_2 \\
		&\text{s.t } x \in C_{\text{DR}}(G_{\pazocal{\hat{P}}_K}), \nonumber 
	\end{aligned}
	\right.
\end{equation}
which, by Definition \ref{DR_core} and considering the data-driven Wasserstein ball as the ambiguity set of each coalition, is equivalent to 
\begin{equation}
	\left\{
	\begin{aligned}
		\min_{x \in \mathbb{R}^N}&	\; ||x||^2_2 \\
		\text{s.t. }& \sum_{i \in \pazocal{N}}x_i = u_\pazocal{N}, \nonumber  \\
		& \sum\limits_{i \in S}x_i \geq \sup_{\mathbb{Q}_S \in \pazocal{\hat{P}}_{K_S}}\mathbb{E}_{\mathbb{Q}_S}[u_S(\xi)],  \ \forall \ S \subset \pazocal{N}.  \nonumber  
	\end{aligned}
	\right.
\end{equation}
Note that calculating the DR core based on the drawn $K_S$ samples from each coalition $S$ boils down to solving  the worst-case expected value problem:
\begin{align}
P'_S:	\sup_{Q_S \in \pazocal{\hat{P}}_{K_S}} \mathbb{E}_{Q_S}[u_S(\xi)], \ \forall \ S \subset \pazocal{N}. \nonumber 
\end{align}
This is an infinite-dimensional optimization problem over probability measures. Under Assumption \ref{properties} each of these programs parameterized by $S$ can be rewritten as a finite-dimensional convex program under Theorem 4.2 in \cite{mohajerin_esfahani_data-driven_2018}, thus taking the dual form:
\begin{equation}
	\left\{
	\begin{aligned}
		\min_{x}&	\; ||x||^2_2 \\
		\text{s.t. }& \sum_{i \in \pazocal{N}}x_i = u_\pazocal{N}\text{ and } \forall \ S \subset \pazocal{N}, \nonumber  \\
		& \sum\limits_{i \in S}x_i  \geq 	\left\{
		\begin{aligned} 
			\min_{\substack{\lambda_S, \ell_{k_S} \\ z_{k_Sm_S}, v_{k_Sm_S}}}& \lambda_S \epsilon_S + \frac{1}{K_S}\sum_{k_S=1}^{K_S}\ell_{k_S} \nonumber  \\
			\text{s.t. }&   [-u_{m_S}]^\ast(z_{k_Sm_S}-v_{k_Sm_S})+ \nonumber \\  
			&\!\!\! \!\!\! \!\!\! \!\!\! \sigma_\Xi(v_{k_Sm_S})- z^\top_{k_S}\xi^{(k_S)} \leq \ell_{k_S} \ \forall k_S, \forall  m_S \nonumber \\
			& \  \ \  \   \|z_{k_Sm_S} \|_* \leq \lambda_S \ \forall k_S, \forall  m_S. \nonumber 
		\end{aligned}
		\right.
	\end{aligned}
	\right.
\end{equation} 
 The problem above can then be rewritten as 
\begin{equation}
	\left\{
	\begin{aligned}
		\min_{x}&	\; ||x||^2_2 \\
		\text{s.t. }& \sum_{i \in \pazocal{N}}x_i = u_\pazocal{N}\text{ and } \forall \ S \subset \pazocal{N}, \nonumber  \\
		&	\left\{
		\begin{aligned} 
			&\exists \ \lambda_S, \ell_{k_S}, z_{k_Sm_S}, v_{k_Sm_S}:  \lambda_S \epsilon_S + \frac{1}{K_S}\sum_{k_S=1}^{K_S}\ell_{k_S} \leq \sum\limits_{i \in S}x_i   \nonumber  \\
			\text{s.t. }&   [-u_{m_S}]^\ast(z_{k_Sm_S}+v_{k_Sm_S})+ \nonumber \\  
			&\sigma_\Xi(v_{k_Sm_S})-z_{k_S}^\top \xi^{(k_S)} \leq \ell_{k_S} \ \forall k_S, \forall  m_S \nonumber \\
			& \  \ \  \   \|z_{k_Sm_S} \|_* \leq \lambda_S \ \forall k_S, \forall  m_S. \nonumber 
		\end{aligned}
		\right.
	\end{aligned} \right.
\end{equation} 
Due to the existential operator we can minimize with respect to $\lambda_S, \ell_{k_S}, z_{k_Sm_S}, v_{k_Sm_S}$ for all $S \subset \pazocal{N}$. This concludes then the proof.  \hfill $\blacksquare$
	
 Note that the additional  decision variables of $P$  correspond to  Wasserstein distance constraints in the primal problems  $P'_S$ for each $S \subset \pazocal{N}$  \cite[Theorem 4.2]{mohajerin_esfahani_data-driven_2018}.

\section{Numerical example}

We consider a stochastic coalitional game with $N=3$ agents and coalitional values $u_{1}= 2+\xi$,   $u_{2}=1.5+\xi$, $u_3=2.5+\xi$ and $u_{12}=6+\xi, u_{23}=6.5+\xi, u_{13}=7+\xi$, where $\xi$ is an uncertain parameter following a Gaussian distribution with mean value $\mu=1$ and variance $\sigma^2=1$, truncated at the interval $\Xi=[0,1]$. We assume that each coalition has a certain level of ambiguity with respect to the distribution of the uncertain parameter. Though our results hold under the consideration of different Wasserstein balls and multi-sample sizes per coalition $S \subset \pazocal{N}$, we assume here the same radius and the same multi-sample size for illustration purposes (see Corollary 1).
 Initially, we consider for each coalition a Wasserstein ball of radius  $\epsilon=0.3$.

\begin{figure}[t]
 	\begin{center}
 		\includegraphics[width=8.5cm,height=5cm]{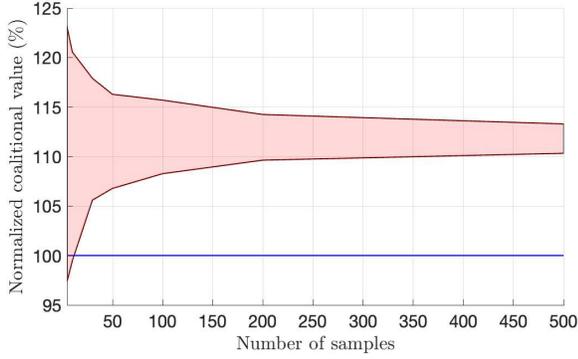}
 	\end{center}
 	\caption{ DR coalitional values for the expectation of $u_S(\xi)$ evaluated over 500 simulations (red shaded area) vs corresponding expected coalitional value  (blue solid line) for $K_S \in \{5, 10, 30, 50, 100, 200, 500\}$ and $\epsilon=0.3$. For $K_S \geq 10$ the DR core is contained with empirical confidence 1 in the expected valued core.}  \label{samples_new}
 \end{figure}

 Figure  \ref{samples_new} focuses on coalition $S=\{1\}$ and shows the range of normalized DR values for the expectation of $u_S(\xi)$,  denoted by $\sup_{\mathbb{Q}_S \in \pazocal{\hat{P}}_{K_S}}\mathbb{E}_{\mathbb{Q_S}}[u_S(\xi)]$ over 500 simulations per multi-sample size $K_S$. We note that a similar behaviour is exhibited in all other coalitions, with the graphs being centered at different values. As the observed pattern is similar across coalitions, it is not shown to avoid repetition. 
 
\begin{figure}[t]
	\begin{center}
		\includegraphics[width=8.5cm,height=5cm]{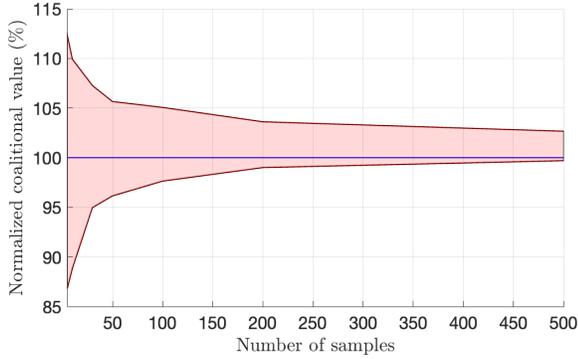}
	\end{center}
	\caption{ DR coalitional values for the expectation of $u_S(\xi)$ evaluated over 500 simulations (red shaded area) vs corresponding expected coalitional value  (blue solid line) for $K_S \in \{5, 10, 30, 50, 100, 200, 500\}$ and $\epsilon=0.03$. It is observed that decrease in $\epsilon$ affects the confidence with which the DR core is contained within the expected value core. At $K_S=500$, however, only a small portion of DR values is below the expected value. }  \label{samples_Sergio}
\end{figure}
 \begin{figure}[t]
 	\centering
 	\includegraphics[width=8.5cm,height=5cm]{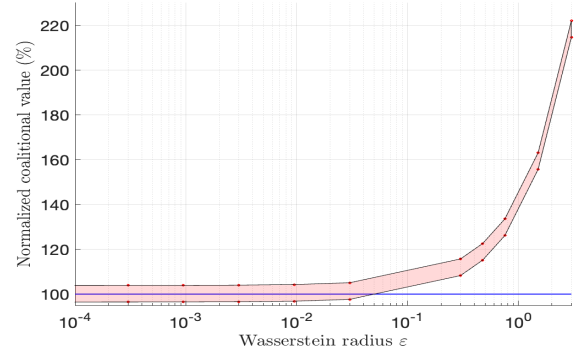}
 	\caption{Coalitional values (red shaded area), among 500 simulations with increasing Wasserstein radius $\epsilon$ for $K_S=100$. As the radius increases, the DR coalitional values of the expectation of $u_S(\xi)$ are above their corresponding expected value (blue solid line) and thus, since the same pattern is repeated across coalitions, the DR core lies within the expected valued core. }  \label{radius}
 \end{figure}

 \begin{figure}[t]
	\centering
	\includegraphics[width=8.5cm,height=5cm]{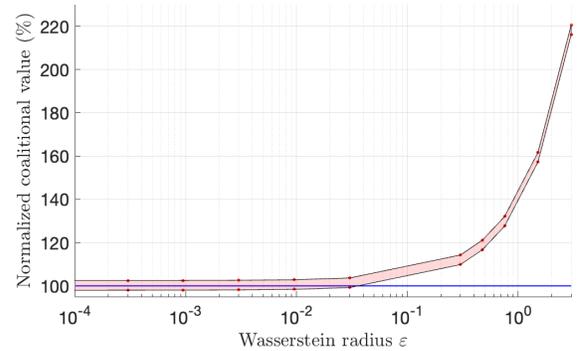}
	\caption{Coalitional values (red shaded area), among 500 simulations with increasing Wasserstein radius $\epsilon$ for $K_S=250$. When samples are drawn from distributions that admit a density with certain concentration properties,  we conjecture that due to the more accurate empirical probability across simulations, the DR coalitional values variability is smaller compared to Figure \ref{radius}. }  \label{radius_Sergio}
\end{figure}

 For the same number of samples per coalition varying in [$5, 500$], we observe that as the number of samples increases, all the DR coalitional values, illustrated by the red shaded area in Figure \ref{samples_new}, are above the corresponding expected value (blue solid line). Since the same pattern is observed across all coalitions, this implies (see Equation  (\ref{prob_guarantess_agent_type}) in the proof of Theorem \ref{main_theorem_dr}) that the DR core is contained within the expected value core  as $K_S$ increases for each $S \subset \pazocal{N}$. This observation is in line with Theorem \ref{DR_guarantees},  since for an increasing multi-sample size per coalition, the confidence $\beta$ in the provided theoretical guarantees tends to 1.  Figure \ref{samples_Sergio} shows that following the same approach as Figure \ref{samples_new} for a significantly smaller Wasserstein radius $\epsilon=0.03$ leads to a smaller empirical confidence that is improved the more samples we obtain.  At $K_S=500$, however, only a small portion of DR values is below the expected value, which implies stability of DR stable allocations in the mean sense with high confidence.

Figure \ref{radius} illustrates the DR coalitional values for the empirical expectation of $u_S(\xi)$ (red shaded areas), over 500 simulations, compared to the corresponding expected value (blue solid line) as the Wasserstein radius $\epsilon$ increases (for a fixed number of samples $K_S=100$). We note again that for a radius $\epsilon$ larger than a certain threshold, all DR coalitional values are above their corresponding expected value and therefore, since this holds for all coalitions when the same number of samples is used, the DR core lies within the expected value core. This behaviour is expected because for a fixed number of samples the larger the Wasserstein ball the more likely it is to include the true probability distribution. As such, obtaining allocations stable in the mean sense can also be achieved, even for a small number of samples, by tuning the Wasserstein radius of each coalition accordingly. We note that the comparison of Figures \ref{samples_new} and \ref{samples_Sergio} is consistent with Figure \ref{radius}.  The percentage of the width of the shaded area acts as an empirical estimate of the confidence. The higher $\epsilon$ or $K_S$ are, the lower $\beta_S$ and as a result the higher the confidence. 

Compared to Figure \ref{radius}, in Figure \ref{radius_Sergio} for each simulation a larger number  $K_S=250$  of samples is generated. Drawing conclusions for the general case is not straightforward, however, when samples are drawn from distributions that admit a density and with certain concentration properties we conjecture that the more samples are used, the more likely it is that the resulting empirical distributions $\hat{P}_{K_S}$ are closer to each other across simulations (i.e., for different multi-samples). As such the centres of the Wasserstein balls would be closer, which in turn implies that the DR value for the expectation of $u_S(\xi)$ would be closer to each other as well. In other words, the higher the number of samples the smaller the variability of the resulting DR value for the expectation of $u_S(\xi)$ across simulations. In line with this intuition, the width of the shaded area in Figure \ref{radius_Sergio} (an empirical estimate of variability) is smaller compared to that of Figure \ref{radius}. \par


\section{Conclusion}
We have introduced the concept of distributionally robust core for coalitional games subject to distributional uncertainty, namely a set of payoff allocations that is robust in the expected value sense. We showed both theoretically and numerically that the concept of distributionally robust stability implies stability in the mean sense as more data becomes available given a certain radius. Furthermore, one can obtain payoff allocations in the expected core by tuning the Wasserstein radius even for a small sample size. 
This paper takes a first step towards studying the class of distributionally robust chance-constrained coalitional games. Future work will  focus on improving the probabilistic guarantees and on designing distributed payoff allocation algorithms.

\bibliographystyle{IEEEtran}
\bibliography{biblio_coop}
\end{document}